\documentclass[11pt]{article}
\usepackage{amsmath,amssymb,latexsym}
\usepackage{theorem}
\newtheorem{theorem}{Theorem}
\newtheorem{lemma}{Lemma}

\newtheorem{corollary}{Corollary}

{\theorembodyfont{\rmfamily} 

\newtheorem{remark}{Remark}} {\theorembodyfont{\slshape} }

\newcommand{\field}[1]{\mathbb{#1}}

\newcommand{\Z}{\field{Z}}
\newcommand{\N}{\field{N}}

\newcommand{\T}{\field{T}}

\renewcommand{\AA}{{\mathcal A}}

\newcommand{\PP}{{\mathcal P}}
\renewcommand{\AA}{{\mathcal A}}

\newcommand{\isdef}{\stackrel{\text{\tiny def}}{=}}

\DeclareRobustCommand{\qed}{%
\ifmmode 
\else \leavevmode\unskip\penalty9999 \hbox{}\nobreak\hfill \fi
\quad\hbox{\qedsymbol}}
\newcommand{\openbox}{\leavevmode
\hbox to.77778em{%
\hfil\vrule
\vbox to.675em{\hrule width.6em\vfil\hrule}%
\vrule\hfil}}
\newcommand{\qedsymbol}{\openbox}
\newcommand{\proofname}{Proof}
\newenvironment{proof}[1][\proofname]{\par
\normalfont \trivlist \item[\hskip\labelsep   \itshape #1. ]
\ignorespaces
}{%
\qed\endtrivlist } 

\def\XXint#1#2#3{{\setbox0=\hbox{$#1{#2#3}{\int}$}
\vcenter{\hbox{$#2#3$}}\kern-.5\wd0}}

\title{On Extensions of a Theorem of Baxter}%
\author{J.\ S.\ Geronimo\\Georgia Institute of Technology, USA \and
 A.\ Mart\'{\i}nez-Finkelshtein\footnote{Corresponding author.
 E-mail \texttt{andrei@ual.es}}\\ Universidad de Almer\'{\i}a, Spain
 }%
\date{\sc Dedicated to Barry Simon on the occasion of his 60th birthday.}%
\begin{document}
\maketitle

\begin{abstract}
We combine the Riemann-Hilbert approach with the techniques of
Banach algebras to obtain an extension of Baxter's Theorem for
polynomials orthogonal on the unit circle. This is accomplished by
using the link between the negative Fourier coefficients of the
scattering function and the coefficients in the recurrence formula
satisfied by these polynomials.
\end{abstract}

\section{Introduction} \label{section:intro}

Recently there has been a upsurge in interest in the theory of
orthogonal polynomials on the unit circle (OPUC). To a large extent
this can be attributed to the two volumes by Barry Simon devoted to
the theory of these polynomials. If $\mu$ is a positive probability
measure supported on the unit circle with an infinite number of
points of increase then a sequence of polynomials
$\varphi_n(z)=\kappa_n z^n + \dots$ of degree $n, n=0,1,\ldots$
satisfying
$$
\frac{1}{2\pi}\, \int_{-\pi}^{\pi}
\varphi_n(e^{i\theta})\overline{\varphi_m(e^{i\theta})}d\mu
=\delta_{n,m},
$$
defines a set of orthonormal polynomials on the unit circle $\T$. We
will assume that the leading coefficient $\kappa_n$ in each
polynomial $\varphi_n$ is positive so the above equation uniquely
specifies the polynomials. It is well know, \cite{Geronimus61},
\cite{Simon05}, \cite{szego:1975} that these polynomials satisfy the
following recurrence formula,
\begin{equation}\label{1.1}
\boldsymbol{\varphi}_{n+1}(z)=A_n(z)\boldsymbol{\varphi}_n(z),    
\end{equation}
with
$$
\boldsymbol{\varphi}_n(z)=\begin{pmatrix}  \varphi_n(z)\\
\varphi_n^*(z)
\end{pmatrix}, \quad
\boldsymbol{\varphi}_0(z)=\begin{pmatrix} 1\\ 1 \end{pmatrix}\,,
$$
and
\begin{equation}\label{1.2}
A_n(z)=\rho_n^{-1}\begin{pmatrix}
 z&-\bar\alpha_n \\ -\alpha_n z&1
\end{pmatrix} .
\end{equation}
In the formula above $\varphi_n^*(z)=z^n\bar \varphi_n(1/z)$ is
called the \emph{reversed} polynomial. The $\alpha_n$ were known as
the \emph{recurrence} or \emph{Verblunsky coefficients}
\cite{Simon05}, and have the property that $|\alpha_n|<1$. Also
$\rho_n=(1-|\alpha_n|^2)^{1/2}=k_n/k_{n+1}$. It is well known that
for any $n\in \N$, $\varphi_n$ has all its zeros strictly inside the
unit circle.

Let us recall (see \cite[Chapter 5]{Simon05}) that a \emph{Beurling
weight} is a two-sided sequence $\nu=\{\nu(n)\}_{-\infty}^{\infty}$
with the properties
\begin{align}
\nu(0)& = 1,\quad \nu(n)\ge 1, \label{1.5}
\\
\nu(n) & =\nu(-n), \label{1.6} \\
\nu(n+m) & \le \nu(n) \nu(m).\label{1.7}
\end{align}
These properties imply the existence of the limit
\begin{equation}\label{1.8}
\lim_{n\to +\infty} \nu(n)^{1/ n}   =R\geq 1\,.
\end{equation}
If $R=1$, then $\nu$ is a \emph{strong} Beurling weight.

Each Beurling weight $\nu$ has associated the Banach spaces
$\ell_\nu$ of two-sided sequences $f=\{f_n \}_{n=-\infty}^{+\infty}$
with
\begin{equation*}\label{1.9}
 \|f\|_\nu   \isdef \sum_{-\infty}^{\infty} \nu(n)|f_n|<\infty;
\end{equation*}
this norm extends naturally to any one-sided sequence by completing
the latter with zeros.

In what follows for a positive measure $\mu$ on $\T$ we will write
\begin{equation}\label{1.4}
d\mu=w\left(\theta\right)d\theta + d\mu_s, 
\end{equation}
after the Lebesgue-Radon-Nikodym decomposition theorem. We will
abuse notation using $w(\theta )$ and $w\left( e^{i \theta}\right)$
interchangeably, expecting that this will not cause confusion.

As noted by B.\ Simon \cite{Simon05}, Baxter's Theorem which relates
the decay rate of the recurrence coefficients $\alpha=\{
\alpha_n\}_{n=0}^{+\infty}$ to the decay rate of the Fourier
coefficients $c=\{ c_n\}_{n=-\infty}^{+\infty}$ of the $\log w$,
\begin{equation}\label{c}
c_n=\frac{1}{2\pi}\, \int_{-\pi}^{\pi} e^{-in\theta } \log w(\theta)
d\theta, \qquad n \in \Z\,,
\end{equation}
plays an important role in the theory:
\begin{theorem}[\cite{Baxter61}]\label{thm1.1}
Suppose $\nu$ is a strong Beurling weight ($R=1$). Then, $\alpha\in
\ell_\nu$ if and only if $\mu_s=0$, $w\in C[-\pi,\pi]$ is strictly
positive, and $c\in \ell_\nu$.
\end{theorem}

In the summer of 2003 Barry Simon asked one of us (JSG) whether it
was possible to extend the above Theorem for the case of Beurling
weights with $R>1$ as had been done for polynomials orthogonal on
the real line (see Geronimo \cite{Geronimo94}). What is needed is to
remove poles from the weight function and to be able to control the
decay of the coefficients of the new weight function. For
polynomials orthogonal on the real line this is accomplished by the
Christoffel-Uvarov formulas (see \cite{szego:1975}, \cite{Uvarov59}
and \cite{Geronimo94}). However such formulas are not available at
present in the case of the unit circle. The necessary control over
the coefficients however can be established using Banach algebra
techniques coupled with those coming from the Riemann-Hilbert
approach (see \cite{Deift05}, \cite{math.CA/0502300}, or for an
alternate approach \cite{Simon05a}). We begin by collecting the
results from scattering theory that will be needed. In particular we
introduce the scattering function $\mathcal S$ and show that its
projection $\PP_-(\mathcal S)$ governs the rate of decay $\alpha$.
We then show that the product of two weights whose coefficients
decay exponentially with the same rate give a weight whose
coefficients decay with the same rate. We then apply this to
Berstein--Szeg\H{o} perturbations of weights. Finally we obtain a
generalization of Baxter's Theorem.

\section{A generalization of Baxter's theorem}

As was noticed by Baxter \cite{Baxter61}, Banach algebras may be
associated with the Beurling weights $\nu$ in the following manner.
Let $a,\ b\in \ell_\nu$; then their convolution is given by
$$
(a*b)(n)=\sum_{k=-\infty}^{\infty}a(k)b(n-k),
$$
which is absolutely convergent by \eqref{1.5} and by \eqref{1.7},
\begin{equation}\label{2.1}
\| a*b \|_\nu=\sum_n \nu(n)\sum_k |a(k)b(n-k)|\le \|a\|_\nu
\|b\|_\nu.
\end{equation}
Thus if we consider the space of functions
\begin{equation}\label{defOfA}
{\mathcal A}_\nu  \isdef  \left\{ f(z)=\sum_{k \in \Z} f_k z^{ k }:
\|f\|_\nu  = \sum_k \nu(k)|f_k|<\infty \right\}\,,
\end{equation}
equation \eqref{2.1} shows that it is closed under multiplication
and so forms an algebra. For $\nu \equiv 1$ we obtain the Wiener
algebra $\mathcal A_1$, containing ${\mathcal A}_\nu$ for any
Beurling weight $\nu$.

On the set of all two-sided sequences $\{d(n)\}_{n\in \Z}$ we define
the projectors $\PP_-$ and $\PP_+$:
$$
\left( \PP_- d\right)(n)=\begin{cases} d(n), & \text{if } n \leq 0,
\\ 0 , & \text{if } n > 0 \,,
\end{cases} \quad \text{and} \quad \left( \PP_+ d\right)(n)=\begin{cases} d(n), & \text{if } n \geq 0,
\\ 0 , & \text{if } n < 0 \,.
\end{cases}
$$
These projectors can be naturally extended to $A_\nu$, and give rise
to two subalgebras associated with ${\cal A}_\nu$: ${\cal A}^+_\nu
\isdef \PP_+(\AA_\nu)$ and ${\cal A}^-_\nu \isdef  \PP_-(\AA_\nu)$.
In other words, $\AA_\nu^{\pm}$ are the set of functions $f\in {\cal
A}_\nu$ whose Fourier coefficients $f_n$ vanish for $n<0$ or $n>0$
respectively. It is easy to see because of \eqref{1.5} that if $f\in
{\cal A}_\nu$ then $f(z)$ is continuous for $1/R\le |z|\le R$ (which
is the maximal ideal space associated to ${\cal A}_\nu$) and
analytic for $1/R< |z|< R$. Likewise if $f$ is in ${\cal A}^+_\nu$
then $f$ is continuous for $|z|\le R$ and analytic for $|z|<R$ (in
which case the series in \eqref{defOfA} is its Maclaurin expansion
convergent at least in this disk), and if $f$ is in ${\cal A}^-_\nu$
then $f$ is continuous for $|z|\ge 1/R$ and analytic for $|z|>1/R$.

If we denote by $\hat \varphi_n(z)\isdef z^{-n} \varphi_n(z)$, then
$\hat \varphi_n \in {\mathcal A}_\nu^-$, and $\varphi_n^* \in
{\mathcal A}_\nu^+$, and we may recast \eqref{1.1} as
\begin{equation}\label{2.2}
G_{n+1}(z)=B_n(z)G_n,
\end{equation}
where
$$
G_n(z)=\begin{pmatrix}   \hat \varphi_n(z)\\
\varphi_n^*(z) \end{pmatrix}, \quad \text{and} \quad
B_n(z)=\rho_n^{-1} \begin{pmatrix} 1&-z^{-(n+1)}\bar\alpha_n\cr
-\alpha_n z^{n+1}&1\end{pmatrix}.
$$

The next lemma is implicitly given by Baxter \cite{Baxter61} and
Simon \cite[Ch.\ 5]{Simon05}. We present its proof for completeness
of the reading.
\begin{lemma}
\label{lemma2.1} Suppose $\nu$ is a Beurling weight, and
$\{\alpha\}_0^{\infty}\in l_\nu$. Then
\begin{align}
\lim_n \varphi_n^* = f_+ & \in {\cal A}_\nu^+, \label{2.4} \\
\lim_n   \hat \varphi_n = f_- & \in {\cal A}_\nu^-,  \label{2.5}
\end{align}
where convergence is understood in the $\|\cdot \|_\nu$ norm. For
$|z|=1$, $f_-(z)=\overline{ f_+(z)}=\bar f_+(1/z)$.

Furthermore, in the decomposition \eqref{1.4}, $d\mu_s=0$, and
\begin{equation}\label{fPlusFMinus}
    w=\frac{1}{f_+ f_-}
\end{equation}
on $\T$; hence, $w$ is a strictly positive continuous function on
$\T$.

\end{lemma}
\begin{proof}
Let us denote by $\Phi_n(z)=\kappa_n^{-1} \varphi_n(z)$ the monic
orthogonal polynomial, and $ \hat\Phi_n(z)=z^{-n} \Phi_n(z)$. Then
from equation \eqref{2.2} we find
\begin{equation}\label{2.6}
    \begin{split}
 \Phi_{n+1}^*(z) & =  \Phi_n^*(z)-\alpha_n\,
 z^{n+1}\hat \Phi_n (z) , \\
 \hat \Phi_{n+1} (z) & =  \hat \Phi_n (z)-\overline{\alpha_n}\, z^{-(n+1)}
 \Phi_n^*(z) .
    \end{split}
\end{equation}
A consequence of \eqref{1.7} is that for $f\in \mathcal A_\nu$  and
$k\in \Z$, $\|z^k f(z)\|_\nu \leq \nu(k)\, \|f\|_\nu$. Thus,
$$
\| \Phi_{n+1}^*\|_\nu \le   \|  \Phi_n^* \|_\nu +
\nu(n+1)|\alpha_n|\, \| \hat \Phi _n \|_\nu .
$$
A similar formula holds for $\hat \Phi_{n+1} $, so that by induction
it follows that
$$
\|g_{n+1}\|_\nu \le \prod_{i=0}^n  (1+\nu(i+1)|\alpha_i|),
$$
with $g_n= \Phi_n^*$ or $\hat \Phi_n$, so that $ \sup_n \|
\Phi_n^*\|_\nu <\infty$, $\sup_n \| \hat \Phi_n\|_\nu < \infty$.


By \eqref{2.6},
\begin{equation*}\label{2.7}
\|  \Phi_n^*-   \Phi_m^*\|_\nu \le \sum_{i=m}^{n-1}
\nu(i+1)|\alpha_i| \|\hat \Phi_i \|_\nu,
\end{equation*}
which yields that $\{\Phi_n^* \}$ is a Cauchy sequence in $\mathcal
A_\nu^+$, and hence has a limit in $\mathcal A_\nu^+$. Since
$$
\kappa_n=\prod_{i=0}^{n-1} \left(1-|\alpha _i|^2 \right)^{-1/2}
$$
converge, \eqref{2.4} follows. Similar reasoning gives \eqref{2.5}.

Since for $|z|=1$, $\hat \varphi_n =\overline{\varphi_n^*}$ we find
for $|z|=1$ that $f_-(z) = \overline{f_+(z)}=\bar f_+(1/z) $.

On the other hand, by \eqref{2.4}--\eqref{2.5},
\begin{equation*}
\lim_n \| |\varphi^*_n|^{2}-  f_- f_+\|_\nu =0.
\end{equation*}
Using  the well know fact from the theory of OPUC
\cite{Geronimus61}, \cite{Simon05}, \cite{szego:1975} that
$\varphi_n^*(z)\ne 0, |z|\le 1$, and for any $h \in C(\T)$, where
$C(\T)$ is the set of continuous functions on the unit circle,
\begin{equation*}\label{2.3}
\frac{1}{2\pi}\, \int_{-\pi}^{\pi} h(\theta)
\frac{d\theta}{|\varphi^*_n(e^{i\theta})|^2} \to \frac{1}{2\pi}\,
 \int_{-\pi}^{\pi} h(\theta) d\mu(\theta)
\end{equation*}
the result follows.
\end{proof}

\begin{remark}
We may recast the last result in terms more familiar to the
orthogonal polynomials community. The \emph{Szeg\H{o} function} of
$w$ (see e.g.\ \cite[Ch.\ X, \S 10.2]{szego:1975}) is defined by
\begin{equation*}\label{standardSzego}
D (z) \isdef \exp\left( \frac{1}{4\pi }\,\int_0^{2\pi} \log w(e^{i
\theta}) \, \frac{e^{i \theta}+z}{e^{i \theta}-z}\, d\theta
\right)\,.
\end{equation*}
This function is piecewise analytic and non-vanishing, defined for
$|z|\neq 1$, and we denote by $D_{\rm i}$ and $D_{\rm e}$ its values
for $|z|<1$ and $|z|>1$, respectively. Then
\begin{equation*}\label{f_and_D}
f_+(z)=D_{\rm i}^{-1}(z)\,, \quad f_-(z)=\ D_{\rm e}(z)\,.
\end{equation*}
Also the last identity of Lemma \ref{lemma2.1} is consistent with
the property
\begin{equation*}\label{symmetry}
\overline{D_{\rm i}\left(\frac{1}{\overline{z}}
\right)}=\frac{1}{D_{\rm e}(z)}\,, \quad |z|>1\,.
\end{equation*}

Lemma \ref{lemma2.1} and equation \eqref{fPlusFMinus} show that $1/w
\in \mathcal{A}_\nu$, which gives a natural meromorphic extension of
$w$ such that $w\neq 0$ for $R^{-1}\leq |z|\leq R$, and this is the
meaning we will be giving to $w$ in what follows.

From the previous works \cite{Deift05}, \cite{Geronimo/Case:79},
\cite{math.CA/0502300}, \cite{Simon05}  it is clear that the
function
\begin{equation}\label{scatteringFunction}
\mathcal{S}=\frac{f_-}{f_+}=D_{\rm i} D_{\rm e}
\end{equation}
plays a prominent role in the theory. By Lemma \ref{lemma2.1} we
have that $|\mathcal S(z)|=1$ for $|z|=1$. Following Faddeyev's
classic article on inverse scattering \cite{Faddeyev63} and Newton's
book on Scattering theory \cite{Newton82} we call the unimodular
function $\mathcal S$ the \emph{scattering function}.

Straightforward computation shows that
$$
\mathcal S(z)=\exp\left( \sum_{k=1}^{+\infty} \left(c_k z^k - c_{-k}
z^{-k} \right)\right)=\exp\left( \sum_{k=1}^{+\infty} \left(c_k z^k
- \overline{c_{k}} z^{-k} \right)\right)\,,
$$
where $c_n$'s are defined in \eqref{c}. Hence, if  $\nu$ is a
Beurling weight, then
\begin{equation}\label{equivalences}
\log w \in \mathcal A_\nu \quad  \Leftrightarrow \quad \log \mathcal
S \in \mathcal A_\nu \quad  \Leftrightarrow \quad  \text{both }
\mathcal S \text{ and } \mathcal S^{- 1} \in \mathcal A_\nu \,.
\end{equation}
In particular, if $\nu$ is a \emph{strong} Beurling weight, then
from Baxter's Theorem it follows that
$$
\log w \in \mathcal A_\nu \quad \Leftrightarrow \quad w \in \mathcal
A_\nu    \quad \Leftrightarrow \quad  \log \mathcal S \in \mathcal
A_\nu \quad \Leftrightarrow \quad \{\alpha\}_0^{\infty}\in l_\nu\,.
$$
\end{remark}

Regarding the case of a Beurling weight with $R>1$, the situation is
different. It was observed by B.\ Simon \cite{Simon05} that the
following three statements are equivalent:
\begin{enumerate}
\item[(i)] the Fourier coefficients of $\log w$ exhibit an exponential
decay;

\item[(ii)] the Fourier coefficients of $w$ exhibit an exponential
decay;

\item[(iii)] $ \alpha _n$ exhibit an exponential
decay.
\end{enumerate}
As \eqref{equivalences} shows, we may replace $w$ by $\mathcal S$ in
(i) and (ii).

However, the decay rates above are not necessary the same; typically
the decay in (ii) and (iii) is faster than in (i). As examples in
\cite[Chapter 7]{Simon05} show, the decay rate in (ii) can be
greater or less than in (iii).

Our goal is to investigate these rates more precisely, extending
Baxter's theorem to the case of $R>1$. Our main results are gathered
in Theorems \ref{thm:crucial}, \ref{thm:product} and
\ref{thm:BaxterGeneral}.

Let us start with the following technical result:
\begin{lemma} \label{lemma:Banach}
Let $h \in \AA_1$, and let $\nu$ be a Beurling weight which
increases on $\Z_+$. Assume that $\PP_-(h) \in \AA_\nu^-$. Then for
a function $g$
$$
g\in \AA_1^+ \text{ or } g \in \AA_\nu^- \quad \Rightarrow \quad
\PP_-(g h) \in \AA_\nu^-\,.
$$
\end{lemma}
\begin{proof}
We begin by proving the first assertion. Set $h=\sum_{k\in \Z} h_k
z^k$ and $g=\sum_{n\ge0} g_n z^n$. Then
$$
\PP_-(gh)=\sum_{i\ge0} d_iz^{-i}\,, \quad \text{with } d_i=
\sum_{n\ge0} g_n h_{-n-i}.
$$
Hence,
$$
\nu(i) |d_i| \leq \nu(i) \sum_{n\ge0} |g_n| \, |h_{-n-i}| \leq
\sum_{n\ge0} \nu(n+i)  |g_n| \, |h_{-n-i}|\,.
$$
Therefore, using Fubini's theorem we get that
$$
\left\| \PP_-(gh)\right \|_\nu \leq \left\| \PP_-(h)\right \|_\nu \,
\left\| g \right \|_1\,.
$$

For the second assertion, write $h_-=\PP_-(h)$ and $h_+=h-h_-\in
\AA_1^+$. Then $\PP_-( g h_-)\in \AA_\nu^-$ by assumption, and
$\PP_-( g h_+)\in \AA_\nu^-$ by the previous argument. Since
$\PP_-(g h)= \PP_-( g h_+)+ \PP_-( g h_-)$, the statement follows.
\end{proof}

\begin{corollary} \label{lemma:Fourier}
Let $g(z)=\sum_{k\in \Z}d_k z^k  \in \mathcal A_1$, and let $\nu$ be
a Beurling weight $\nu$ which increases on $\Z_+$. For a rational
function $r$ denote $r(z) g(z)=\sum_{k\in \Z}\hat d_k z^k$. If $r$
has no poles in $1/R\leq r\leq 1$ then
$$
\sum_{k < 0} d_k z^k \in \mathcal A_\nu^- \quad \Rightarrow \quad
\sum_{k < 0} \hat d_k z^k \in \mathcal A_\nu^-\,.
$$
\end{corollary}

Let us recall some facts that we will need further (see
\cite{math.CA/0502300} for details). If
$$
\rho ^{-1} \isdef \varlimsup_{n\to \infty} |\alpha _n|^{1/n}<1\,,
$$
then (see \cite{Nevai/Totik:89}) also
\begin{equation*}
\rho^{-1} =  \inf\{0<r<1:\, f_-(z) \text{ is holomorphic in } |z|>r
\}\,.
\end{equation*}
For the scattering function $\mathcal S$ defined in
\eqref{scatteringFunction}, let
$$
\mathcal S(z)=\sum_{k\in \Z} d_k z^k
$$
be its Laurent expansion. An important fact, established
independently in \cite{Deift05}, \cite{Simon05a} and
\cite{math.CA/0502300} (cf.\ Corollary 2 therein) is that for $n\in
\N$,
\begin{equation}\label{formulaForAlphas}
\overline{\alpha _n}=-d_{-n-1} + e_n\,, \quad \text{with } |e_n|\leq
C\, r^{3n}\,,
\end{equation}
for an arbitrary $r$ such that $1/\rho <r\leq 1$, and where $C>0$
does not depend on $r$ or $n$.

\begin{remark}
Since $\mathcal S(z)|=1$ for $|z|=1$, we may state this result
equivalently in terms of the coefficients of $\PP_+(1/\mathcal S)$,
see \cite{Martinez-Finkelshtein05}.
\end{remark}

The above formular indicate the close connection between the rate of
decay of $\alpha _n$'s and the decay of the negative Fourier
coefficients of the scattering function $\mathcal S$, which is made
explicit in:
\begin{theorem} \label{thm:crucial}
Let $d\mu=w\left(\theta\right)d\theta$ with $w\in \AA_1$ and $\nu$
be a Beurling weight. Then
$$
\alpha \in \ell_\nu \quad \Leftrightarrow \quad \PP_-(\mathcal{S})
\in \AA_\nu^-\,.
$$
\end{theorem}
\begin{proof}
Since for any Beurling weight, $\AA_\nu \subset \AA_1$, Baxter's
theorem (Theorem \ref{thm1.1}) shows that $f_\pm$ exist, $f_\pm \in
\AA_1 ^\pm$, and therefore $\mathcal S \in \AA_1$. Furthermore, for
any strong Beurling weight the results follow from Theorem
\ref{thm1.1} and Lemma \ref{lemma:Banach}, so we consider only the
case when $R>1$.

Assume that $\alpha \in \ell_\nu$, then $1<R\leq \rho $. For an
sufficiently small $\varepsilon >0$ let us take $r=(1+\varepsilon
)/R \in (1/\rho , 1)$. Since we may take $r$ arbitrary close to
$1/R$, assumption $\alpha \in \ell_\nu$ and \eqref{formulaForAlphas}
imply that $\{ d_k \}_{k<0} \in \ell_\nu$, so that $\PP_-(\mathcal
S)\in \mathcal A_\nu^-$.

For converse, if $\PP_-(\mathcal S)\in \mathcal A_\nu^-$, then Lemma
\ref{lemma:Banach} implies that $ f_-\in \mathcal A_\nu^-$, showing
that again $\rho >1$. Now an application of formula
\eqref{formulaForAlphas} gives the assertion.
\end{proof}

A consequence of the previous theorem is the following result which
shows that the product of two measures whose coefficients are in
$\mathcal A_\nu$ is a new measure whose coefficients are also in
$\mathcal A_\nu$.
\begin{theorem} \label{thm:product}
For $d\mu_i=w_i\left(\theta\right)d\theta$, $i=1, 2, 3$, with
$w_i\in \AA_1$ denote by $\alpha ^{(i)}$ the corresponding sequences
of the recurrence coefficients. Assume that $w_3=w_1 w_2$, and let
$\nu$ be a Beurling weight. Then
$$
\alpha ^{(1)}, \alpha ^{(2)} \in \ell_\nu \quad \Rightarrow \quad
\alpha ^{(3)}\in \ell_\nu \,.
$$
\end{theorem}
\begin{proof}
Let us denote by $\mathcal S_i$ the scattering function
corresponding to $\mu_i$. By Baxter's theorem, $\mathcal S_1,
\mathcal S_2 \in \AA_1$, and by Theorem \ref{thm:crucial}, for $i=1,
2$,
$$
\alpha ^{(i)} \in \ell_\nu \quad \Rightarrow \quad \mathcal S_i^-
\isdef  \PP_-(\mathcal S_i) \in \AA_\nu^-, \quad \mathcal S_i^+
\isdef \mathcal S_i- \mathcal S_i^- \in \AA_1^+.
$$
Hence,
$$
\PP_-(\mathcal S_3)=\PP_-(\mathcal S_1 \mathcal S_2)=\PP_-(\mathcal
S_1^- \mathcal S_2^- + \mathcal S_1^+ \mathcal S_2^-+ \mathcal S_1^-
\mathcal S_2^+) \in \AA_\nu^-,
$$
where we have used Lemma \ref{lemma:Banach}. The result now follows
from Theorem \ref{thm:crucial}.
\end{proof}

Now we are in the position of discussing several corollaries of the
previous results.

\begin{corollary} \label{thm:deift}
Let $\nu$ be a Beurling weight, and for the decomposition
\eqref{1.4}, $d\mu_s=0$ and $w\in \mathcal A_\nu$. If $w(z)\neq 0$
in $R^{-1}\leq |z|\leq R$, where $R$ is defined in \eqref{1.8}, then
$\{\alpha\}_0^{\infty}\in \ell_\nu$.
\end{corollary}
This result for complex weights with zero winding numbers has been
proved in \cite{Deift05}.
\begin{proof}
According to \eqref{equivalences} and Theorem \ref{thm:crucial},
under assumptions of the corollary,
$$
\log w \in \mathcal A_\nu \quad  \Rightarrow \quad \log \mathcal S
\in \mathcal A_\nu \quad  \Rightarrow \quad   \mathcal S  \in
\mathcal A_\nu \quad  \Rightarrow \quad \PP_-(\mathcal S ) \in
\mathcal A_\nu \quad  \Rightarrow \quad  \alpha \in \ell_\nu \,.
$$
\end{proof}

Another straightforward consequence of Theorem \ref{thm:product} is
that the multiplication of the original weight by a
Bernstein-Szeg\H{o} weight does not affect the rate of decay of the
recurrence coefficients:
\begin{corollary} \label{corollary:Bernstein}
Let $d\mu=w\left(\theta\right)d\theta$,  with $w\in \AA_1$, and let
$p$ be a polynomial with no zeros on $\T$. Denote by $\alpha$ and
$\hat \alpha$ the sequences of the recurrence coefficients
corresponding to measures $d\mu= w(\theta)d\theta$ and $d\hat \mu=
|p|^{-2}w(\theta)d\theta$, respectively. Then for any Beurling
weight $\nu$,
$$
\alpha \in \ell_\nu \quad \Rightarrow \quad \hat \alpha \in
\ell_\nu\,.
$$
\end{corollary}
The proof is an immediate consequence of Theorem \ref{thm:product},
by observing that for the weight $|p|^{-2}$ the corresponding
recurrence coefficients vanish for a sufficiently large $n$, and
hence belong to $\ell_\nu$ for any Beurling weight $\nu$.

Finally we formulate an extension of Baxter's theorem (Theorem
\ref{thm1.1}) for the case $R>1$.
\begin{theorem}\label{thm:BaxterGeneral}
Suppose that $\nu$ is a Beurling weight with $R>1$.

If there exists a polynomial $p(z)\neq 0$ for $|z|=1$, such that for
\begin{equation}\label{def_wHat}
\hat w (z)\isdef |p(z) |^2 w(z) \,, \quad |z|=1 \,,
\end{equation}
we have $\log (\hat w ) \in \mathcal A_{\nu}$, then
$\{\alpha\}_0^{\infty}\in \ell_\nu$.

Conversely, if $\{\alpha\}_0^{\infty}\in \ell_\nu$ and $1/w \neq 0$
for $|z|=R$ then there exists a polynomial $p$  such that for $\hat
w (z)$ defined by \eqref{def_wHat} we have $\log (\hat w ) \in
\mathcal A_{\nu}$.
\end{theorem}
\begin{proof}
Let us prove the first implication. Assume that $\log (\hat w ) \in
\mathcal A_{\nu}$. By Corollary \ref{thm:deift}, if $\hat \alpha$
are the recurrence coefficients corresponding to $\hat w$, then
$\{\hat \alpha\}_0^{\infty}\in \ell_\nu$. But $w=|p|^{-2} \hat w$,
and it remains to apply Corollary \ref{corollary:Bernstein}.

Conversely, let $\{\alpha\}_0^{\infty}\in \ell_\nu$. Since by Lemma
\ref{lemma2.1}, $1/w =f_+ f_-\in \mathcal A_\nu$, and $1/w(z)\neq 0$
for $|z|=R$, then $f_+$ has a finite number of zeros, $\zeta_1,
\dots, \zeta_m$, in $1<|z| <R$ (denoted with account of their
multiplicity). Let
$$
p(z)=\prod_{k=1}^m (z-\zeta_k)\,.
$$
Then $\hat w$ satisfies the conditions of Corollary \ref{thm:deift},
which concludes the proof.
\end{proof}

\section*{Acknowledgement}

One of the authors (J.S.G.) would like to thank Barry Simon for
suggesting the extension of Baxter's Theorem to him.

The research of J.S.G. was supported, in part, by a grant from the
National Science Foundation. The research of A.M.F.\ was supported,
in part, by a research grant from the Ministry of Science and
Technology (MCYT) of Spain, project code BFM2001-3878-C02, by Junta
de Andaluc\'{\i}a, Grupo de Investigaci\'{o}n FQM229, and by Research Network
on Constructive Complex Approximation (NeCCA), INTAS 03-51-6637.

Both A.M.F.\ and J.S.G.\ acknowledge also a partial support of NATO
Collaborative Linkage Grant ``Orthogonal Polynomials: Theory,
Applications and Generalizations'', ref. PST.CLG.979738.

%

\end{document}